\documentclass{article}

\usepackage[english]{babel}

\usepackage[a4paper]{geometry}

\usepackage{amsmath}
\usepackage{amssymb}
\usepackage{amsthm}
\usepackage{graphicx}
\usepackage[colorlinks=true, allcolors=blue]{hyperref}
\usepackage{subfigure}
\usepackage{color}

\newtheorem{remark}{Remark}

\title{On the computation of the SVD of \\Fourier submatrices}
\author{Simon Dirckx, Daan Huybrechs and Robbe Ongenae}

\begin{document}
\maketitle

\begin{abstract}
Contiguous submatrices of the Fourier matrix are known to be ill-conditioned. In a recent paper in SIAM Review A. Barnett has provided new bounds on the rate of ill-conditioning of the discrete Fourier submatrices~\cite{barnett2022fouriersubmatrix}. In this paper we focus on the corresponding singular value decomposition. The singular vectors go by the name of \emph{periodic discrete prolate spheroidal sequences} (P-DPSS). The singular values exhibit an initial plateau, which depends on the dimensions of the submatrix, after which they decay rapidly. The latter regime is known as the \emph{plunge region} and it is compatible with the submatrices being ill-conditioned. The discrete prolate sequences have received much less study than their continuous counterparts, prolate spheroidal wave functions, associated with continuous Fourier transforms and widely studied following the work of Slepian in the 1970's. In this paper we collect and expand known results on the stable numerical computation of the singular values and vectors of Fourier submatrices. We illustrate the computations and point out a few applications in which Fourier submatrices arise.
\end{abstract}

\section{The Discrete Fourier Transform}

A common definition of the discrete Fourier transform of a vector $\{ x_k \}_{k=0}^{N-1}$ of length $N$ is
\begin{equation*}
 X_k = \sum_{n=0}^{N-1} x_n e^{-2 \pi i k n/ N}, \qquad k = 0,\ldots,N-1.
\end{equation*}
Its inverse is
\begin{equation*}
 x_n = \frac{1}{N} \sum_{k=0}^{N-1} X_k e^{2 \pi i k n/ N}, \qquad n = 0,\ldots,N-1.
\end{equation*}
Other definitions may differ in the position of the $1/N$ term, and in the sign of the exponent in the complex exponential.

For convenience of implementation we prefer to work with vectors $\mathbf{x}$ and $\mathbf{X}$ in $\mathbb{C}^N$ with indices ranging from $1$ to $N$, as is standard in linear algebra, although many of the following formulas are mathematically more elegant with indices starting at $0$. With this convention, the $N \times N$ Fourier matrix $F \in \mathbb{C}^{N \times N}$ is given by
\begin{equation*}
F_{j,k} = e^{-2 \pi i (j-1) (k-1) / N} = \omega^{-(j-1)(k-1)}, \qquad 1 \leq j,k \leq N,
\end{equation*}
in which we have defined the \emph{twiddle factor}
\begin{equation*}
    \omega = e^{2 \pi i / N}.
\end{equation*}
The inverse DFT amounts to a matrix-vector product with $F^*/N$, i.e., we have
\begin{equation*}
\mathbf{X} = F \mathbf{x} \quad \textrm{and} \quad \mathbf{x} = \frac{1}{N} F^* \mathbf{X}.
\end{equation*}

\subsection{Submatrices of the Fourier matrix}

We focus on the contiguous submatrix $A \in \mathbb{C}^{p \times q}$ starting from the $(1,1$) entry of $F$,
\begin{equation*}
 A_{j,k} = F_{j,k}, \qquad 1 \leq j \leq p, \quad 1 \leq k \leq q,
\end{equation*}
with $1 \leq p,q \leq N$. Any other contiguous submatrix of $F$ with the same dimensions as $A$ is easily obtained by diagonal scalings. Shifting the submatrix along rows and columns corresponds to multiplication by diagonal matrices to the right and to the left  respectively.

Let us be precise. We define the diagonal matrix $D_n \in \mathbb{C}^n$, with diagonal entries
\[
 [D_n]_{j,j} = \omega^{-(j-1)}, \qquad j=1,\ldots,n.
\]
Consider the Fourier submatrix $B$ of dimension $p \times q$ but starting from index $(j_0,k_0)$ of $F$, and taking into account periodicity if the submatrix `wraps around',
\[
 B_{j,k} = F_{1+\textrm{mod}(j_0+j-2,N),1+\textrm{mod}(k_0+k-2,N)},\qquad 1 \leq j \leq p, \quad 1 \leq k \leq q.
\]
Then we have
\begin{equation}\label{eq:B}
 B = \omega^{-(j_0-1)(k_0-1)} \, D_p^{k_0-1} \, A \, D_q^{j_0-1}
\end{equation}
and, conversely,
\[
A = \omega^{(j_0-1)(k_0-1)} \, D_p^{-(k_0-1)} \, B \, D_q^{-(j0-1)}.
\]
Hence, in the following, without loss of generality we consider only $A$.

\subsection{The singular value decomposition of $\boldmath{A}$}

Our goal is the computation of the singular value decomposition of $A$,
\begin{equation}\label{eq:svd}
    A = U \Sigma V^*.
\end{equation}
With $r = \min(p,q)$, we aim to compute the reduced SVD with $U \in \mathbb{C}^{p \times r}$, $\Sigma \in \mathbb{C}^{r \times r}$ and $V \in \mathbb{C}^{q \times r}$. A full SVD can also be computed by generating $|p-q|$ additional orthogonal vectors in $U$, if $p > q$, or in $V$, if $q > p$. A typical profile of singular values is illustrated further on in Fig.~\ref{fig:spectra}(a) for $N=128$ and a few choices of $p$ and $q$.

\section{Computation of the SVD}

\subsection{Historical perspective}
We collect and generalize the formulas for Periodic Discrete Prolate Spheroidal Sequences (P-DPSS) from~\cite{grunbaum1981toeplitz,jain1981extrapolation,xu1984periodicprolate}. The P-DPSS were introduced by Gr{\"u}nbaum~\cite{grunbaum1981toeplitz}, and Jain and Ranganath~\cite{jain1981extrapolation}. Earlier studies were devoted to the (continuous) prolate spheroidal wave functions~\cite{slepian1978prolate} and discrete (but non-periodic) prolate spheroidal sequences~\cite{slepian1978prolateV}. Both of these are eigenfunctions of a bandlimiting operator, but are better computed as eigenfunctions of a related commuting operator: see~\cite{osipov2013prolate} for a comprehensive treatment. To that end, for the periodic discrete setting of P-DPSS,  Gr{\"u}nbaum~\cite{grunbaum1981toeplitz} derived a tridiagonal matrix that commutes with a Toeplitz matrix, associated with the Fourier submatrix $A$. Xu and Chamzas pointed out that two singular values of $A$, $0$ and $\sqrt{N}$, may have higher multiplicity~\cite{xu1984periodicprolate}. This means that the corresponding singular vectors are not uniquely defined. We follow the methodology of~\cite{xu1984periodicprolate} and define singular vectors of $A$ using the unique (up to normalization) eigenvectors of the tridiagonal matrix. Asymptotic and non-asymptotic properties of the P-DPSS were studied further in~\cite{edelman1999futurefft,zhu2018eigenvalue,barnett2022fouriersubmatrix}.

Gr{\"u}nbaum studied Fourier series of odd length, i.e., with symmetric frequencies, and correspondingly assumed both $p$ and $q$ to be odd. A similar assumption is made in the explicit formulas in~\cite{matthysen2016fastfe} in the context of the Fourier extension approximation scheme. Xu and Chamzas adopt the assumption that one of $p$ or $q$ is odd~\cite{xu1984periodicprolate}. Jain and Ranganath in~\cite[\S IX.B]{jain1981extrapolation} do not make any such assumptions in the definition of P-DPSS, but they do not supply formulas. Edelman et al. in~\cite{edelman1999futurefft} study square submatrices and assume that $p$ divides $N$. Zhu et al.~\cite{zhu2018eigenvalue} and Barnett~\cite{barnett2022fouriersubmatrix} make no assumptions on $N$, $p$ or $q$, but they do not supply computational formulas. Moreover, due to differences in notation and indexing (starting from $0$ or symmetric around $0$) in all these references, implementation of the formulas require some attention and care.

The formulas in this section are a generalization of those given in the cited references \cite{grunbaum1981toeplitz,matthysen2016fastfe}, to accommodate any combination of $N$, $p$, and $q$. They are provided for completeness and ease of implementation. They do not contain new insights and are simply stated without proof. The main complication is the correct treatment of even $p$ and $q$.

\subsection{The periodic prolate matrix}\label{ss:sinc}

The left and right singular vectors of $A$ are eigenvectors of $A^* A$ and $A A^*$ respectively. These are complex matrices. The computations simplify by considering a more symmetric matrix $C$ first, which is such that $C^* C$ and $C C^*$ are real.

Matrix $C \in \mathbb{C}^{p \times q}$ is similar to the shifted submatrix $B$ given by~\eqref{eq:B}, though with possibly fractional shifts if either $p$ or $q$ is even:\footnote{We have to fix branch cuts in case $p$ or $q$ is even. The $(k,k)$ entry of $D_p^{-\frac{(q-1)}{2}}$ here is $\omega^{(k-1)(q-1)/2}$ and it is best implemented that way explicitly. Numerically raising the diagonal matrix $D_p$ to the power $-\frac{(q-1)}{2}$ may yield a different outcome, due to choosing a different square root. A similar comment applies to powers of $D_q$.}
\begin{equation}\label{eq:C}
 C = \omega^{-\frac{(p-1)}{2}\frac{(q-1)}{2}} \, D_p^{-\frac{(q-1)}{2}} \, A \, D_q^{-\frac{(p-1)}{2}}.
\end{equation}
Equivalently, its entries are
\[
 C_{j,k} = \omega^{-\left(j-\frac{p+1}{2}\right)\left(k-\frac{q+1}{2}\right)}, \qquad 1 \leq j \leq p, 1 \leq k \leq q,
\]
and we have the inverse relation
\begin{equation}\label{eq:C_to_A}
 A = \omega^{\frac{(p-1)}{2}\frac{(q-1)}{2}} \, D_p^{\frac{(q-1)}{2}} \, C \, D_q^{\frac{(p-1)}{2}}.
\end{equation}

The real Hermitian matrices $C^*C$ and $C C^*$ are both \emph{prolate matrices}. The \emph{periodic prolate matrix} $S(p,q) \in \mathbb{R}^{q \times q}$ is defined by
\begin{equation}\label{eq:sinc}
 [S(p,q)]_{j,k} = \frac{\sin\left(\frac{p(j-k) \pi}{N}\right)}{\sin\left(\frac{(j-k) \pi}{N}\right)}, \qquad 1 \leq j,k \leq q, j \neq k,
\end{equation}
for the off-diagonal entries, and (by H\^opital's rule)
\[
 [S(p,q)]_{j,j} = p
\]
along the diagonal. The matrix $S(p,q)$ can also be viewed as a \emph{Dirichlet kernel matrix}\footnote{For each $n$ the Dirichlet kernel $\mathcal{D}_n$ is a radial kernel given (up to a factor) by $\mathcal{D}_n(r)=\frac{\sin{((2n+1)r/2)}}{\sin(r/2)}$.}. We have
\begin{equation}\label{eq:C_sinc}
    C^* C = S(p,q) \quad \mbox{and} \quad C C^* = S(q,p).
\end{equation}
Hence, possibly up to a normalization factor, the left and right singular vectors of $C$ are given by the eigenvectors of $S(q,p)$ and $S(p,q)$, respectively.

In the special case where $q=N$, one can verify that $S(p,N)$ is an orthogonal projection onto $\text{span}\{(\omega^{-(j-1)((k-1)/2)})_{j=1}^{N}|k\in \{2,4,\ldots,p\}\cup\{2N-p+2,2N-p+4,\ldots,2N\}\}$. However, in general, $S(p,q)$ is not associated with an orthogonal projection.

\begin{figure}[t]
    \centering
    \subfigure[Singular values of $A$]{
    \includegraphics[width=6cm]{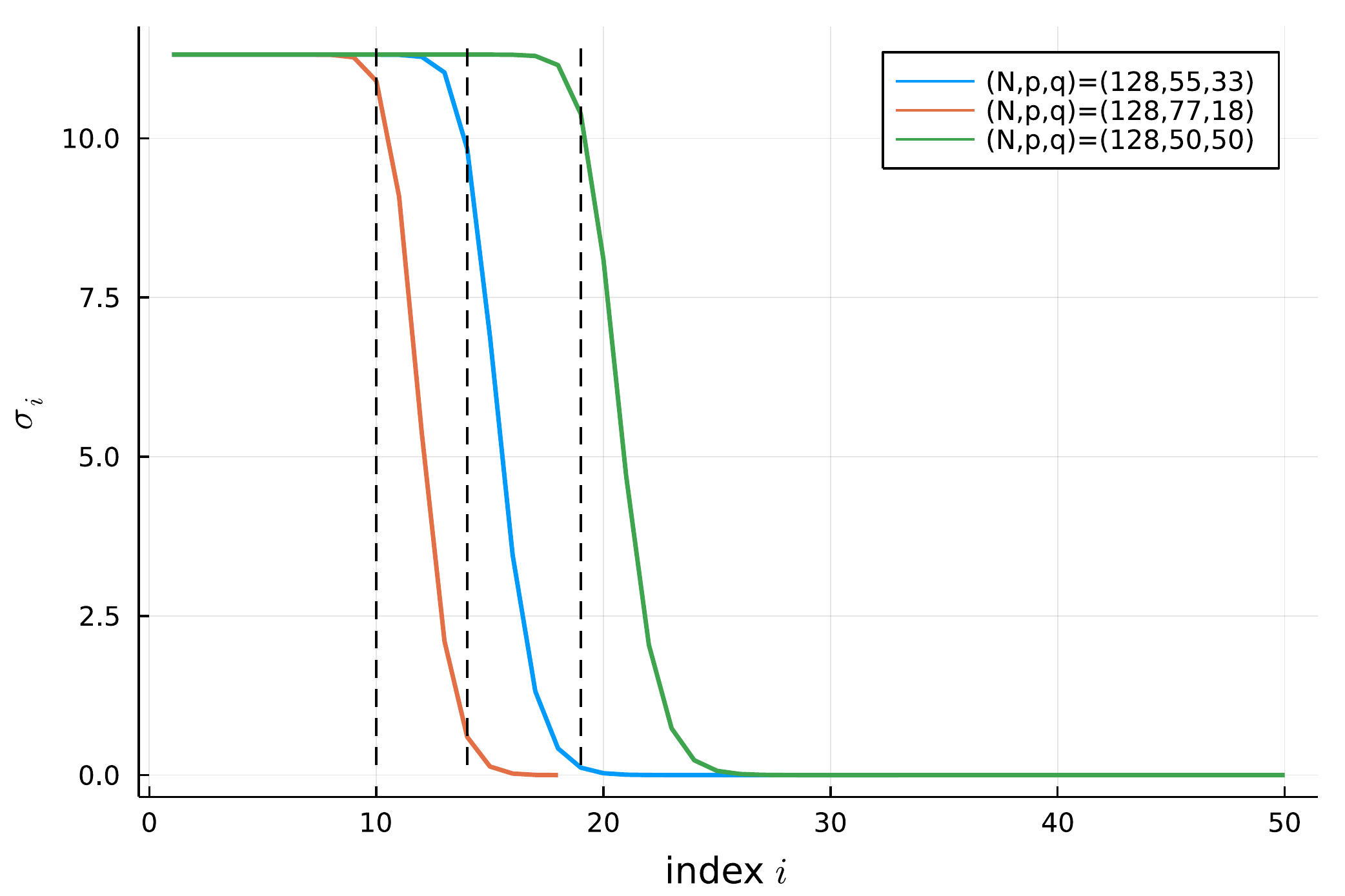}
    }
    \subfigure[Eigenvalues of $J(p,q)$]{
    \includegraphics[width=6cm]{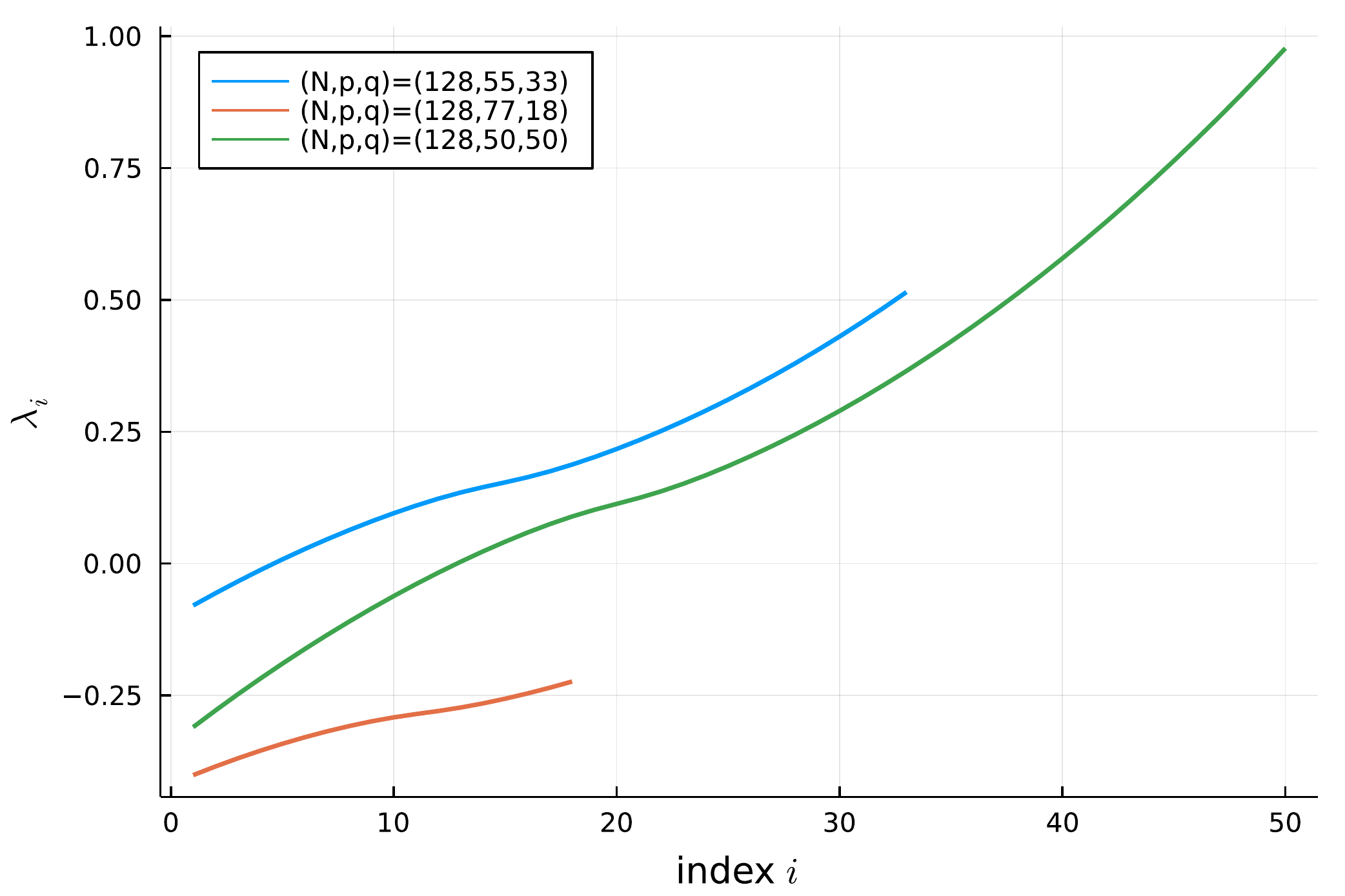}
    }
    \subfigure[Eigenvalues of $J(q,p)$]{
    \includegraphics[width=6cm]{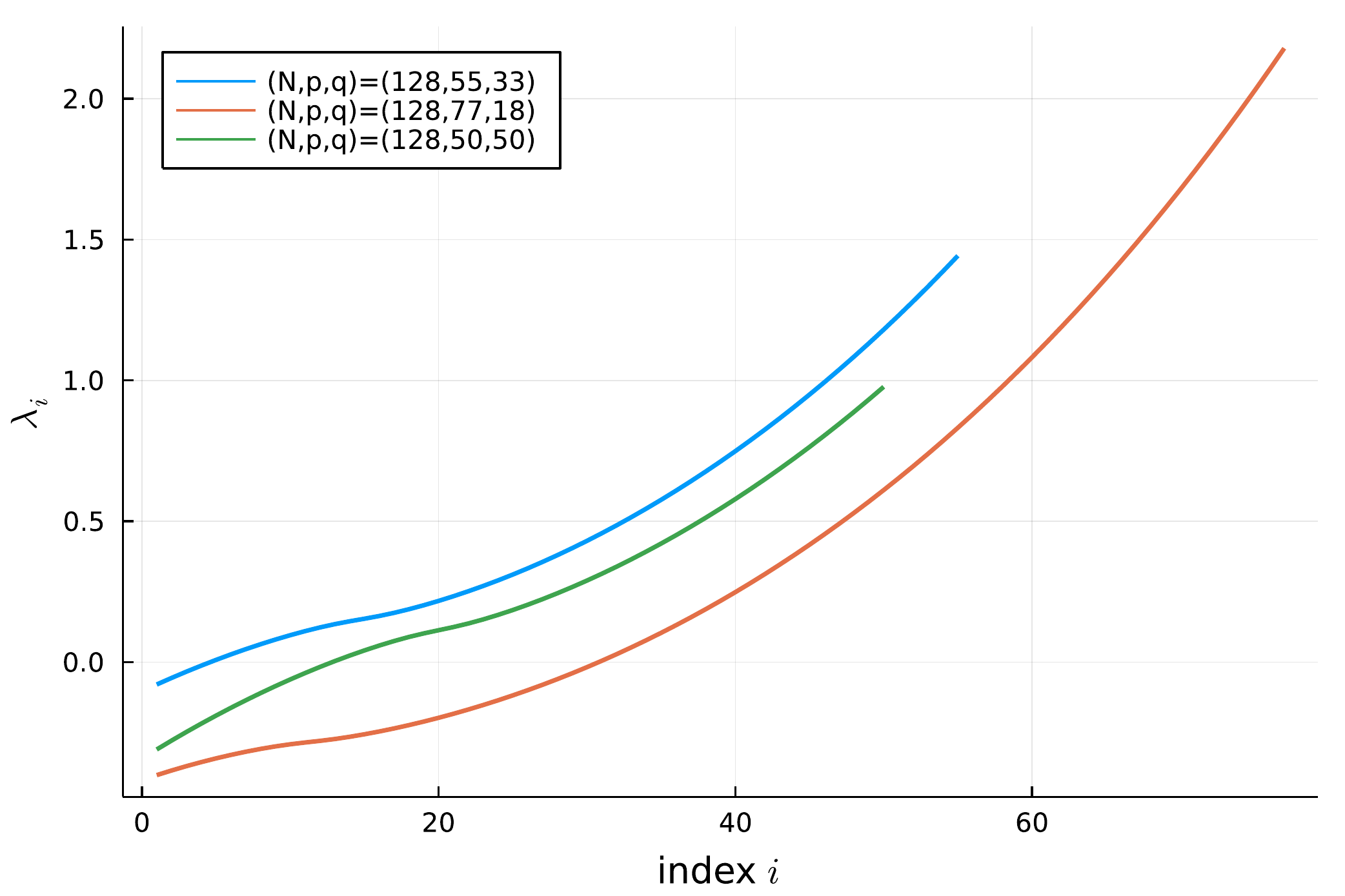}
    }
    \caption{The singular values of Fourier submatrix $A$ cluster towards $\sqrt{N}$ and towards $0$ (top-left panel). The plateau contains approximately $pq/N$ singular values, shown by dashed lines. In contrast, the spectra of the symmetric tridiagonal matrices $J(p,q)$ and $J(q,p)$ are free of clusters and contain only simple eigenvalues. All eigenvectors can be computed to high accuracy.}
    \label{fig:spectra}
\end{figure}

\subsection{Two commuting tridiagonal matrices}\label{ss:tridiag}

The spectrum of the periodic prolate matrix is clustered. This renders the computation of its eigenvectors ill-conditioned, a problem that also appears in the continuous case. The enabling observation for computational purposes, first established in this context by Gr{\"u}nbaum~\cite{grunbaum1981toeplitz}, is that $S(p,q)$ commutes with a real and symmetric tridiagonal matrix $J(p,q) \in \mathbb{R}^{q \times q}$.

The diagonal entries of this tridiagonal matrix are
\[
    [J(p,q)]_{k,k} = \cos\left(\frac{\pi (2k-q-1)}{N}\right)  \, \cos\left(\frac{p \pi}{N}\right), \qquad  1 \leq k \leq q.
\]
The first diagonal above the main diagonal is given by
\[
    [J(p,q)]_{k,k+1} = -\sin\left( \frac{\pi k}{N}\right)  \sin\left(\frac{\pi(q-k)}{N}\right), \qquad 1 \leq k \leq q-1.
\]
The subdiagonal elements $[J(p,q)]_{k+1,k} = [J(p,q)]_{k,k+1}$ are defined by symmetry. Since these elements are all nonzero, the \emph{Sturm sequence property} (see~\cite[Theorem 8.5.1]{golub1996matrixcomputations}) ensures that $J(p,q)$ has no eigenvalues of multiplicity greater than one.

Similarly, we can consider $J(q,p$). We have the commutation relations
\[
 J(p,q)S(p,q) - S(p,q)J(p,q) = 0 \quad \mbox{and} \quad J(q,p)S(q,p) - S(q,p)J(q,p) = 0.
\]
This, combined with the simple spectrum of $J(p,q)$, implies that the eigenvectors of $J(p,q)$ are also eigenvectors of $S(p,q)$. The converse might not be true, since $S$ may have repeated eigenvalues. The spectrum of $J(p,q)$ is also free of clusters, and hence the eigenvalue problem of $J(p,q)$ is better conditioned than that of $S(p,q)$~\cite{xu1984periodicprolate}. The above also holds if we flip $p$ and $q$.

\subsection{The singular vectors of $\boldmath{A}$}

Let $\tilde{v}_k \in \mathbb{R}^q$ be an eigenvector of $J(p,q)$ with eigenvalue $\tilde{\lambda}_k \in \mathbb{R}$,
\begin{equation}\label{eq:eigen_Jpq}
 J(p,q) \tilde{v}_k = \tilde{\lambda}_k \tilde{v}_k.
\end{equation}
All quantities are real-valued owing to the symmetry of $J(p,q)$. Since $J(p,q)$ and $S(p,q)$ commute, we also have
\[
S(p,q) \tilde{v}_k = \lambda_k \tilde{v}_k
\]
for some $\lambda_k \in \mathbb{R}$ which in general differs from $\tilde{\lambda}_k$, but which is still real. Consequently, by~\eqref{eq:C_sinc},
\begin{equation}\label{eq:C_star_C}
C^* C \tilde{v}_k = \lambda_k \tilde{v}_k.
\end{equation}
We can conclude that the corresponding singular value of $C$ is $\sigma_k = \sqrt{\lambda_k}$. However, $\sigma_k$ is not computed this way, since we solve the eigenvalue problem~\eqref{eq:eigen_Jpq} for the tridiagonal matrix rather than~\eqref{eq:C_star_C}. There is no known relation between $\tilde{\lambda}_k$ and $\lambda_k$.

Moving to $A$, from~\eqref{eq:C} we find that
\[
 A^* A D_q^{-\frac{(p-1)}{2}} \tilde{v}_k = \lambda_k D_q^{-\frac{(p-1)}{2}} \tilde{v}_k.
\]
Hence, up to normalization, we find that
\begin{equation}\label{eq:vk}
 v_k = D_q^{-\frac{(p-1)}{2}} \tilde{v}_k
\end{equation}
is a right singular vector of $A$.

For the left singular vectors we proceed similarly. Thus, $\tilde{u}_k \in \mathbb{R}^p$ is an eigenvector of $J(q,p)$ with eigenvalue $\tilde{\lambda}_k \in \mathbb{R}$ (different from the eigenvalue in~\eqref{eq:eigen_Jpq} above),
\begin{equation}\label{eq:eigen_Jqp}
 J(q,p) \tilde{u}_k = \tilde{\lambda}_k \tilde{u}_k.
\end{equation}
The corresponding eigenvalue for $S(q,p)$ is $\lambda_k \in \mathbb{R}$. We have
\begin{equation*}
C C^* \tilde{u}_k = \lambda_k \tilde{u}_k.
\end{equation*}
The corresponding singular value is $\sigma_k = \sqrt{\lambda_k}$. From~\eqref{eq:C} we find that
\[
 A A^* D_p^{\frac{(q-1)}{2}} \tilde{u}_k = \lambda_k D_p^{\frac{(q-1)}{2}} \tilde{u}_k.
\]
Hence, again up to normalization, the left singular vector $u_k$ of $A$ is
\begin{equation}\label{eq:uk}
 u_k = D_p^{\frac{(q-1)}{2}} \tilde{u}_k.
\end{equation}

\subsection{The singular values of $\boldmath{A}$}

It remains to determine the singular value and a suitable normalization for $u_k$ and $v_k$. This is achieved most easily by computing the projections
\begin{equation}\label{eq:sigma_projection}
 \hat{\sigma}_k = u_k^* A v_k, \qquad 1 \leq k \leq \min(p,q)
\end{equation}
where we have ordered the eigenvectors $v_k, u_k$ by decreasing corresponding eigenvalues. The computed value $\hat{\sigma}_k$ may be complex, yet $\sigma_k$ should be real by the definition of the SVD. Thus, the singular value is
\[
 \sigma_k = | \hat{\sigma}_k|.
\]
We can incorporate the phase difference into one of the singular vectors, and in our implementation we arbitrarily choose to do so in $u_k$:
\begin{equation}\label{eq:normalization}
 u_k \to \frac{|\hat{\sigma}_k|}{\hat{\sigma}_k} u_k, \qquad 1 \leq k \leq \min(p,q).
\end{equation}

\subsection{Efficient matrix-vector product using an FFT of length $\boldmath{N}$}
\label{ss:fft}

The computation of the singular values via~\eqref{eq:sigma_projection} requires $\min(p,q)$ matrix vector products with $A$, which is a costly affair. However, since $A$ is a submatrix of the DFT matrix of dimension $N$, an efficient matrix-vector product is also available using the FFT of length $N$~\cite{matthysen2016fastfe}. The multiplication $A v$ can be computed by zero-padding $v$ to a vector of length $N$ (note that $v$ has length $q$), followed by the FFT, followed by restricting to the first $p$ entries. We write that as
\begin{equation}\label{eq:RFE}
    A = R_p F E_q,
\end{equation}
with $E_q \in \mathbb{R}^{N \times q}$, $F$ the DFT matrix of dimension $N$, and $R_p \in \mathbb{R}^{p \times N}$.

Moreover, we have the relation
\begin{equation}\label{eq:from_v_to_u}
 A v_k = \sigma_k u_k, \qquad 1 \leq k \leq \min(p,q).
\end{equation}
Hence, once $v_k$ is determined, both $\sigma_k \in \mathbb{R}$ and $u_k$ can be determined from~\eqref{eq:from_v_to_u} using the FFT. The normalization of $u_k$ is fixed by requiring $\sigma_k$ to be real. This is an alternative to~\eqref{eq:eigen_Jqp}--\eqref{eq:normalization}.

Whether the FFT-based approach is more efficient than computing the matrix-vector product with $A$ directly depends on the relative sizes of $p$ and $q$ versus $N$: the matrix vector product has cost ${\mathcal O}(pq)$, the FFT has cost ${\mathcal O}(N \log N)$. In some applications $N$ may be much larger than $p$ and $q$, for example in the parallel FFT computation of~\cite{edelman1999futurefft}.

\subsection{Computational complexity}

Summarizing, we solve~\eqref{eq:eigen_Jpq} for $\tilde{v}_k$, and use~\eqref{eq:vk} to find the right singular vectors of $A$. At this point we have a choice. We can solve~\eqref{eq:eigen_Jqp} for $\tilde{u}_k$, use~\eqref{eq:uk} to find the left singular vectors of $A$, compute the singular value from~\eqref{eq:sigma_projection} and finally normalize $u_k$ following~\eqref{eq:normalization}. The alternative to~\eqref{eq:eigen_Jqp}--\eqref{eq:normalization} is to use~\eqref{eq:from_v_to_u} to find both $\sigma_k$ and $u_k$ simultaneously.

The computational complexities of these steps are as follows:
\begin{itemize}
 \item The eigenvalue problem~\eqref{eq:eigen_Jpq} for the symmetric tridiagonal matrix  $J(p,q)$ has complexity ${\mathcal O}(q^2)$, and the corresponding problem~\eqref{eq:eigen_Jqp} for $J(q,p)$ is ${\mathcal O}(p^2)$.
 \item The diagonal scaling in~\eqref{eq:vk} is just ${\mathcal O}(q)$, and that of~\eqref{eq:uk} is ${\mathcal O}(p)$.
 \item Determining the singular values from~\eqref{eq:sigma_projection} or~\eqref{eq:from_v_to_u} requires $\min(p,q)$ matrix-vector products with $A$. The cost is:
 \begin{itemize}
     \item ${\mathcal O}(\min(p,q) pq)$ when multiplying with the dense matrix $A$, and
     \item ${\mathcal O}(\min(p,q) N \log N)$ when using the FFT for the matrix-vector products with $A$.
 \end{itemize}
\end{itemize}
The computation of the singular values is the most expensive step in all cases and it determines the overal computational complexity.

\begin{remark}
These formules have been stated for the reduced SVD. The full SVD of $A$ can also be computed by computing $|p-q|$ additional orthogonal vectors from the full eigenvalue decomposition of $J(p,q)$ if $q > p$ or from that of $J(q,p)$ if $p > q$. This does not change the complexity, as there remain at most $\min(p,q)$ nonzero singular values.
\end{remark}

\begin{remark}\label{rem:plunge}
The methodology can also be adapted to compute a part of the spectrum, for example the singular values and vectors associated with the plunge region.
We refer to~\cite[\S4.1]{matthysen2016fastfe} for the details of this computation and for an estimate of the index range that corresponds to the plunge region. It is a topic of considerable interest that the size of the plunge region scales logarithmically with $N$~\cite[Theorem 1]{edelman1999futurefft}.
\end{remark}

\section{Illustrations and applications}

\begin{figure}
    \centering
    \includegraphics[width=8cm]{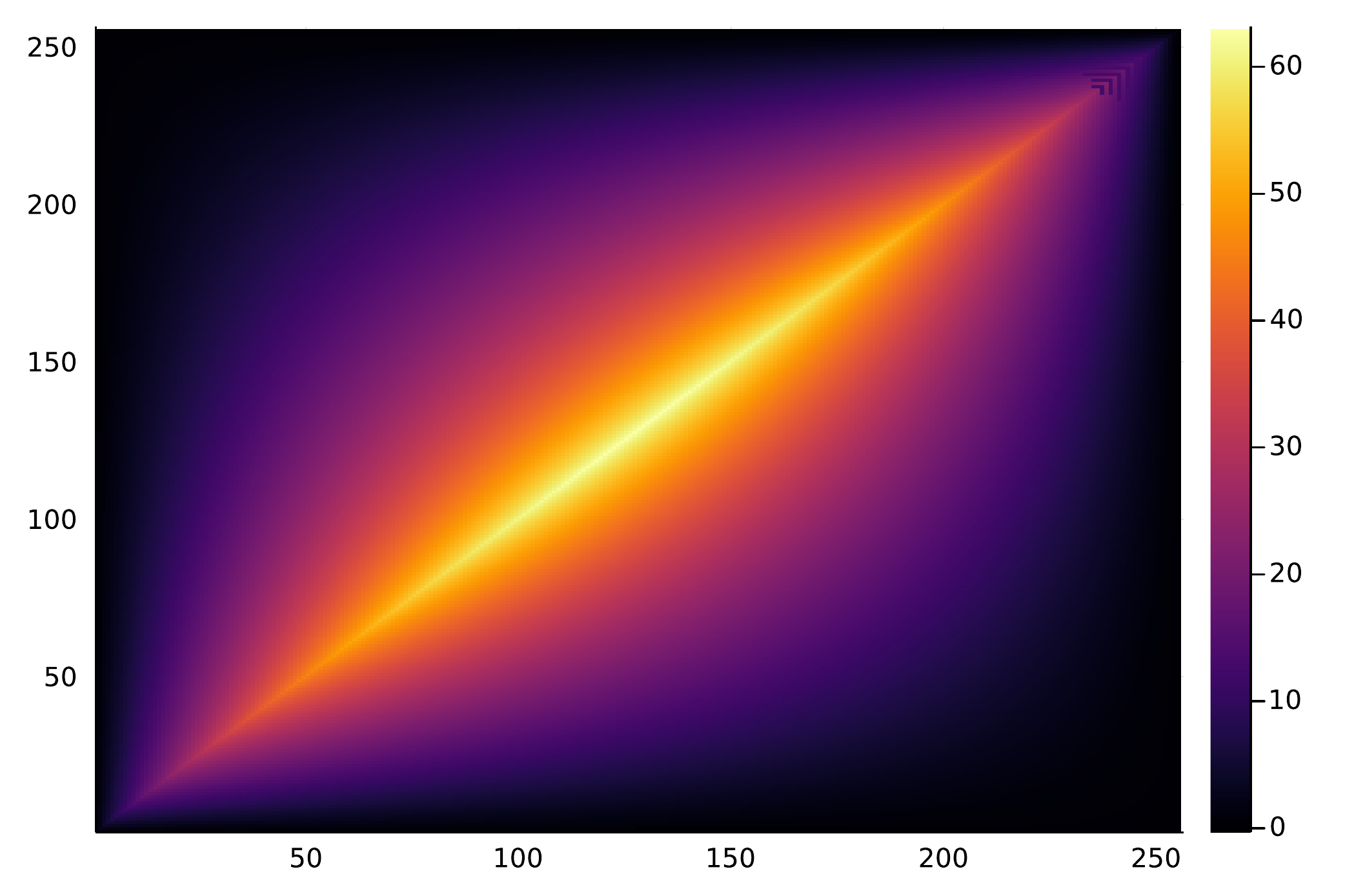}
    \caption{Condition number of all Fourier submatrices of the DFT matrix of dimension 256, shown in base-10 logarithmic scale, with $p$ and $q$ varying between $1$ and $256$. After Barnett~\cite{barnett2022fouriersubmatrix}.}
    \label{fig:cond_N256}
\end{figure}

For the specific application in the Fourier extension approximation scheme we refer to~\cite{matthysen2016fastfe,matthysen2017fastfe2d}, and for a generic algorithm for solving linear systems with DFT submatrices to~\cite{coppe2020az}. Here, we illustrate the properties and review a few applications which are more directly related to the SVD of $A$.

\subsection{The spectra of $\boldmath{S(p,q)}$ and $\boldmath{J(p,q)}$}

Some examples of the spectra of the matrices involved are shown in Fig.~\ref{fig:spectra} for the case $N=128$. The singular values of $A$ in panel (a) exhibit the familiar plateau and plunge regions, varying with $p$ and $q$. The eigenvalues of the tridiagonal matrices $J(p,q)$ and $J(q,p)$ are shown in panels (b) and (c), respectively. The values themselves are not relevant here as they are not used: what matters is that they are simple and, in comparison to the singular values of $A$, well-separated. 

\subsection{Exponential ill-conditioning of $\boldmath{A}$}

For the computation of the condition number of $A$, we have only to compute the smallest and largest singular values. In turn, we compute the smallest and largest eigenvalues of $J(p,q)$, making use of the fact that it is a tridiagonal matrix. This can be achieved with a computational cost proportional to $\max(p,q)$. This enables the computation of the condition number in linear time. This compares favourably to a computation based on the full SVD of $A$, which requires ${\mathcal O}(pq^2)$ or ${\mathcal O}(p^2q)$ operations.

We do note that efficiency is not the main computational problem here. Instead it is the exponential ill-conditioning of $A$. The methods of this paper enable the accurate computation of all the singular vectors, but not of all the singular values, in particular obviously not the ones smaller than machine precision. Hence, computing condition numbers larger than $1/\varepsilon$, where $\epsilon$ is floating point precision, remains elusive unless higher-precision arithmetic is used.

The condition numbers of all submatrices of the DFT matrix of length $N=256$ are shown in Fig.~\ref{fig:cond_N256}. Similar figures were shown in~\cite{barnett2022fouriersubmatrix} for $N=8$, $N=16$ and $N=32$. In our experiment we have first ran all computations in standard double precision. All condition numbers estimated to be greater than $10^{13}$ have been recomputed using higher precision arithmetic with about $70$ digits. Owing to the efficient algorithm this remains quite cheap. The limitation to $N=256$ is not due to computational cost, but due to the fact that the smallest singular values approach $10^{-70}$ in this example, underflowing even in higher precision calculations. In the figure, the largest condition number is $1.5\times 10^{63}$ with $p=q=128$ (and $N=256$).

\subsection{Low-rank relation between DFTs of different length}

\begin{figure}[t]
    \centering
    \subfigure[$N=128$]{
    \includegraphics[width=6cm]{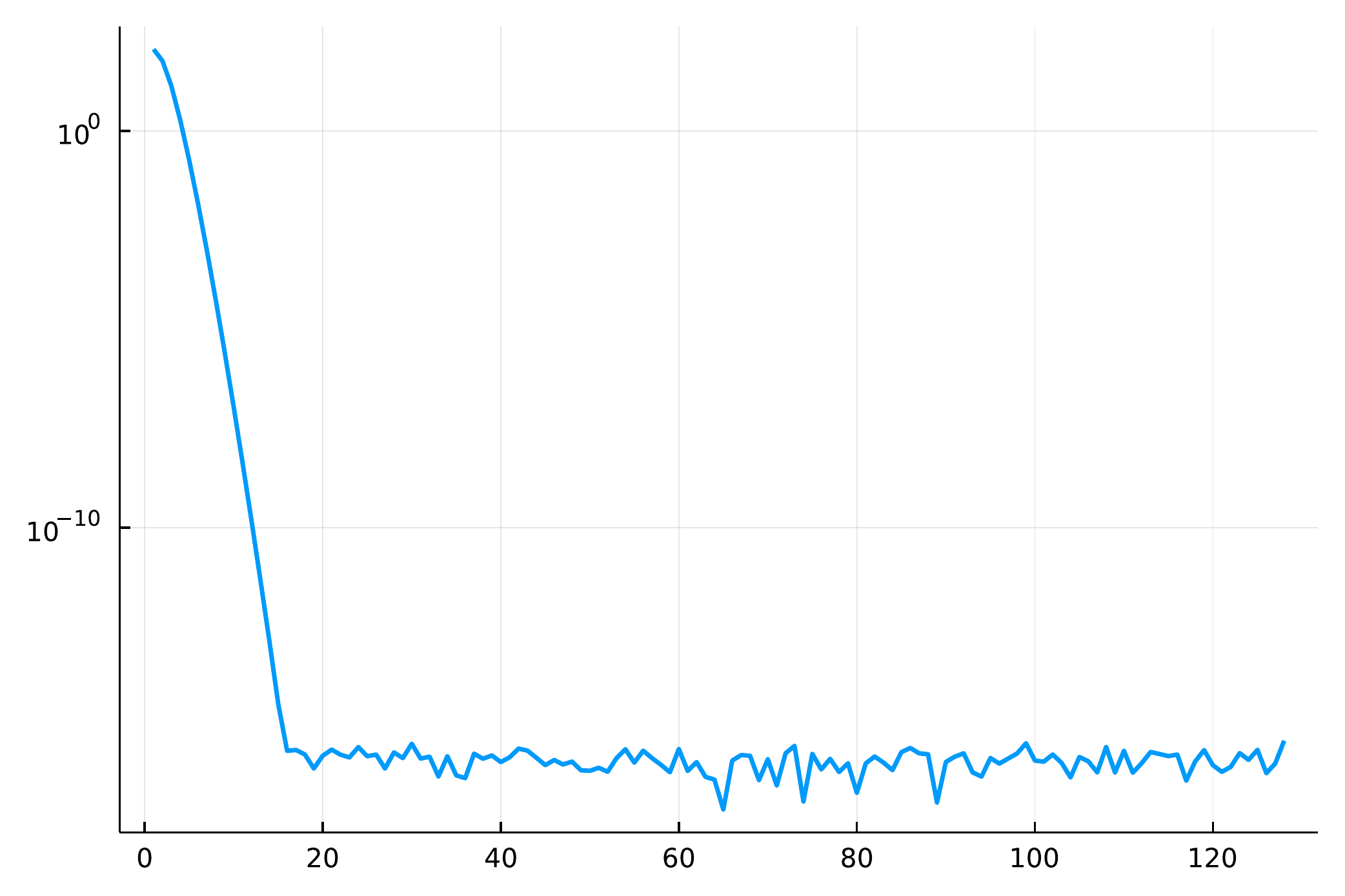}
    }
    \subfigure[$N=512$]{
    \includegraphics[width=6cm]{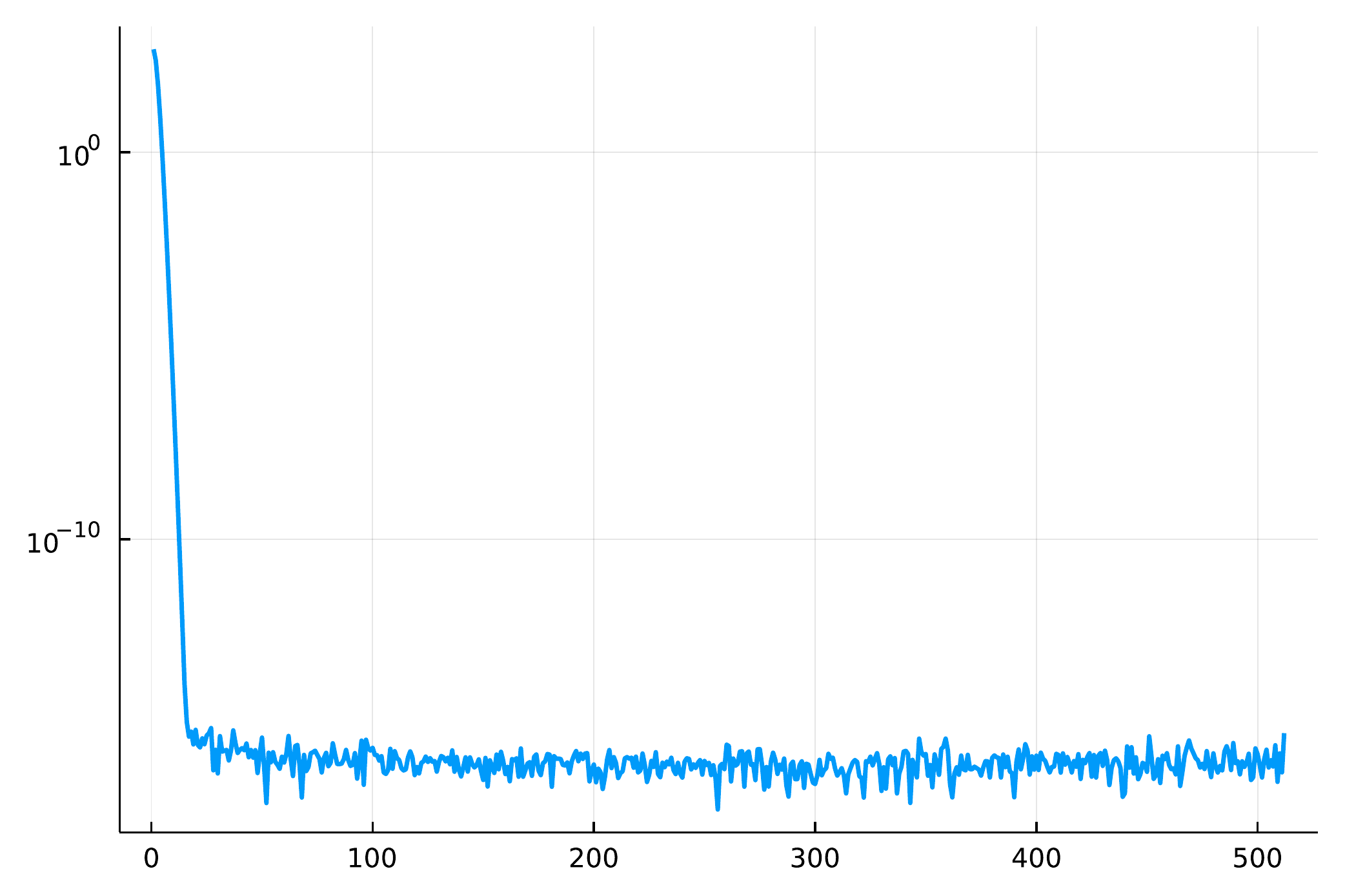}
    }
    \caption{The singular values of $H$ in~\eqref{eq:hadamard} for different values of $N$. The matrices have no plateau, only a plunge region, which scales logarithmically with $N$: $H$ is low rank.}
    \label{fig:spectra_H}
\end{figure}

The efficient evaluation of a DFT is usually based on recursive computation of DFTs of shorter length, the celebrated example being the FFT itself. The recursive relations are typically based on number-theoretic properties of $N$ in combination with algebraic properties of the transform.

Using low-rank Fourier submatrices, we can also establish an approximate numerical link between DFTs of different length. This is based on a Hadamard product, the element-wise product of matrices, which was used for the computation of non-uniform FFTs in~\cite{townsend2018nonuniform}. We consider just one example, namely between the DFT of length $N$ and of length $N+1$.

Consider the $N\times N$ submatrix $G$ of the DFT matrix of dimension $N+1$, and compare to the DFT matrix $F$ of dimension $N$. Thus,
\[
 F_{j,k} = e^{-2 \pi i \frac{(j-1) (k-1)}{N}} \quad \mbox{and} \quad G_{j,k} = e^{-2 \pi i \frac{(j-1) (k-1)}{N+1}}.
\]
Using the Hadamard product notation $\circ$, we write the link as
\begin{equation}\label{eq:hadamard}
 F = G  \circ H.
\end{equation}
Elementwise, this corresponds to $F_{j,k} = G_{j,k}H_{j,k}$. From the explicit formulas of the matrices $G$ and $F$ above, we deduce that 
\[
 H_{j,k} = e^{-2 \pi i \frac{(j-1) (k-1)}{N(N+1)}}.
\]
This is essentially due to the identity
\[
 \frac{1}{N} = \frac{1}{N+1} + \frac{1}{N (N+1)}.
\]
Hence, $H$ is precisely the $N\times N$ submatrix of the much larger DFT matrix of dimension $N(N+1)$.

The singular values of $H$ exhibit no plateau, since the general formula $\frac{pq}{\dim(F)}$ here becomes, with $p=q=N$ and $\dim(F) = N(N+1)$,
\[
 \frac{N^2}{N(N+1)} = \frac{N}{N+1} < 1.
\]
There is only a plunge region with ${\mathcal O}(\log N)$ values larger than any given fixed threshold.

From~\eqref{eq:hadamard} we can also write
\[
 G = F \circ \overline{H}.
\]
The main computational technique of~\cite{townsend2018nonuniform} is to expedite matrix-vector products of the form $(A \circ B)\mathbf{x}$. This is shown to be possible at least if $A$ has a fast matrix-vector product and $B$ has low rank. If $B$ has rank $r$, then the matrix-vector product with $A\circ B$ requires $r$ matrix-vector products with $A$. That is exactly the setting of this subsection. However, this observation is mainly of mathematical interest. It is unlikely to result in any computational advantage compared to existing FFT methods, because even with $r = {\mathcal O}(\log N)$, the Hadamard matrix-vector product remains expensive in comparison.

\subsection{Localization properties of P-DPSS}

\begin{figure}[t]
    \centering
    \includegraphics[width=\linewidth]{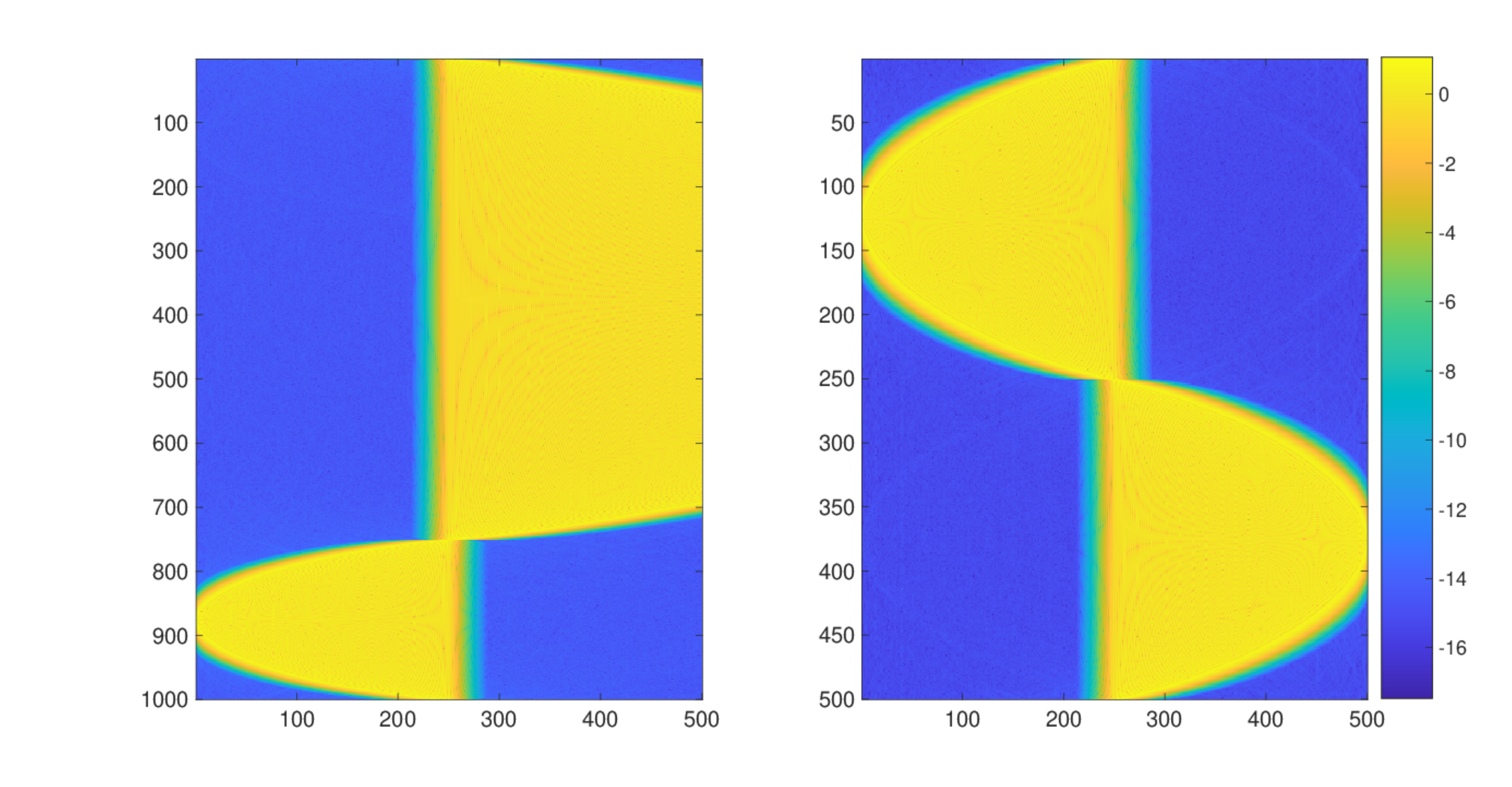}
    \caption{Illustration of the separation property of the prolate sequence matrices $U$ and $V$ in the frequency domain. Plotted are magnitudes $\log|F_pU(i,j)|$, and $\log|F_qV(i,j)|$ over the indices in the matrices $F_pU$ and $F_qV$ respectively. The figure shows that $F_pU$ and $F_qV$ each consist of three regions: purely low frequency, purely high-frequency and a transitory regime in-between. Here $N=2000$, $p=1000$, $q=500$.}
    \label{fig:separation}
\end{figure}

Finally, we illustrate the time-frequency properties of the singular vectors, i.e., of the periodic discrete prolate spheroidal sequences. There are three regimes, associated with the three regimes of the singular values shown in Fig.~\ref{fig:spectra}(a). The singular vectors corresponding to the larger singular values are approximately localized in the lower part of the spectrum, while those corresponding to small singular values have high-frequency content. The plunge region is associated with a transitory regime covering the full spectrum.

These properties are conveniently visualized simply by taking the column-wise FFT of the matrices $U$ and $V$ of length $p$ and $q$ respectively. This is shown in Fig.~\ref{fig:separation}. The separation between low and high frequency is the ratio $p/N = 1/2$ in the left panel, and $q/N=1/4$ in the right panel. In both panels, the left part of the figure corresponds to large singular values and has the exact opposite localization properties as the right part of the figure corresponding to small singular values. The plunge region appears as a small band in the middle covering the full frequency spectrum.

\section*{Acknowledgements}

The authors are grateful to Alex Barnett and Nick Trefethen for their suggestions an discussions on the topic, and to Emiel Dhondt for his contributions in the implementation.

\bibliographystyle{abbrv}
\bibliography{refs}

\begin{thebibliography}{10}

\bibitem{barnett2022fouriersubmatrix}
A.~H. Barnett.
\newblock How exponentially ill-conditioned are contiguous submatrices of the
  {F}ourier matrix?
\newblock {\em SIAM Review}, 64:105--131, 2022.

\bibitem{coppe2020az}
V.~Copp{\'e}, D.~Huybrechs, R.~Matthysen, and M.~Webb.
\newblock The {AZ} algorithm for least squares problems with a known incomplete
  generalized inverse.
\newblock {\em SIAM J. Mat. Anal. Appl.}, 46:1237--1259, 2020.

\bibitem{edelman1999futurefft}
A.~Edelman, P.~McCorquodale, and S.~Toledo.
\newblock The future {F}ast {F}ourier {T}ransform?
\newblock {\em SIAM J. Sci. Comput.}, 20(3):1094--1114, 1999.

\bibitem{golub1996matrixcomputations}
G.~H. Golub and C.~F. van Loan.
\newblock {\em Matrix computations}.
\newblock Johns Hopkins University Press, Baltimore, MD, 3rd edition, 1996.

\bibitem{grunbaum1981toeplitz}
F.~A. Gr{\"u}nbaum.
\newblock Eigenvectors of a {T}oeplitz matrix: discrete version of the prolate
  spheroidal wave functions.
\newblock {\em SIAM J. Alg. Disc. Meth.}, 2(2):136--141, 1981.

\bibitem{jain1981extrapolation}
A.~Jain and S.~Ranganath.
\newblock Extrapolation algorithms for discrete signals with application in
  spectral estimation.
\newblock {\em IEEE Trans. Acoust., Speech, Signal Processing},
  ASSP-29(4):830--845, 1981.

\bibitem{matthysen2016fastfe}
R.~Matthysen and D.~Huybrechs.
\newblock Fast algorithms for the computation of {F}ourier extensions of
  arbitrary length.
\newblock {\em SIAM J. Sci. Comput.}, 38(2):A899--A922, 2016.

\bibitem{matthysen2017fastfe2d}
R.~Matthysen and D.~Huybrechs.
\newblock Function approximation on arbitrary domains using fourier frames.
\newblock {\em SIAM J. Numer. Anal.}, 56:1360--1385, 2018.

\bibitem{osipov2013prolate}
A.~Osipov, V.~Rokhlin, and H.~Xiao.
\newblock {\em Prolate spheroidal wave functions of order zero}.
\newblock Springer, New York, 2013.

\bibitem{townsend2018nonuniform}
D.~Ruiz-Antol{\'i}n and A.~Townsend.
\newblock A nonuniform {F}ast {F}ourier {T}ransform based on low rank
  approximation.
\newblock {\em SIAM J. Sci. Comput.}, 40(1):A529--A547, 2018.

\bibitem{slepian1978prolateV}
D.~Slepian.
\newblock Prolate spheroidal wave functions, {F}ourier analysis, and
  uncertainty {V}: the discrete case.
\newblock {\em The Bell System Technical Journal}, 57:1371--1430, 1978.

\bibitem{slepian1978prolate}
D.~Slepian and H.~Pollak.
\newblock Prolate spheroidal wave functions, {F}ourier analysis, and
  uncertainty—{I}.
\newblock {\em Bell Systems Tech. J.}, 40(1):43--63, 1961.

\bibitem{xu1984periodicprolate}
W.~Y. Xu and C.~Chamzas.
\newblock On the periodic discrete prolate spheroidal sequences.
\newblock {\em SIAM J. Appl. Math.}, 44(6):1210--1217, 1984.

\bibitem{zhu2018eigenvalue}
Z.~Zhu, S.~Karnik, M.~A. Davenport, J.~Romberg, and M.~B. Wakin.
\newblock The eigenvalue distribution of discrete periodic time-frequency
  limiting operators.
\newblock {\em IEEE Signal Processing Letters}, 25(1):95--99, 2018.

\end{thebibliography}

\end{document}